\documentclass[12pt]{amsart}
\usepackage[a4paper,asymmetric]{geometry}
\usepackage{latexsym}
\usepackage{amsthm}
\usepackage{amssymb}
\usepackage{amsfonts}
\usepackage{amsmath}
\newtheorem{theorem}{Theorem}[section]
\newtheorem{thm}[theorem]{Theorem}

\newtheorem{lem}[theorem]{Lemma}
\newtheorem{proposition}[theorem]{Proposition}

\newtheorem{corollary}[theorem]{Corollary}
\newtheorem{assumption}[theorem]{Assumption}
\theoremstyle{definition}

\newtheorem{defn}[theorem]{Definition}
\newtheorem{ex}[theorem]{Example}
\theoremstyle{remark}
\newtheorem{remark}[theorem]{Remark}
\newtheorem{rem}[theorem]{Remark}
\numberwithin{equation}{section}

 \DeclareMathAlphabet{\mathpzc}{OT1}{pzc}{m}{it}
 
 \newcommand{\E}{\mathbb{E}}            
 
 \newcommand{\e}{\varepsilon}
 \newcommand{\p}{\partial}

 \newcommand{\Ll}{\langle}
 \newcommand{\Rr}{\rangle}

 \newcommand{\N}{\mathbb{N}}
 \newcommand{\R}{\mathbb{R}}
 
 \newcommand{\PP}{\mathbb{P}}
 \newcommand{\mcl}{\mathcal}
 
 \newcommand{\Be}{\begin{equation}}
 \newcommand{\Ee}{\end{equation}}
 \newcommand{\Bs}{\begin{split}}
 \newcommand{\Es}{\end{split}}
  \newcommand{\Bes}{\begin{equation*}}
 \newcommand{\Ees}{\end{equation*}}
 \newcommand{\BT}{\begin{thm}}
 \newcommand{\ET}{\end{thm}}
 \newcommand{\Bp}{\begin{proof}}
 \newcommand{\Ep}{\end{proof}}
 \newcommand{\BL}{\begin{lem}}
 \newcommand{\EL}{\end{lem}}
 \newcommand{\BP}{\begin{proposition}}
 \newcommand{\EP}{\end{proposition}}
 \newcommand{\BC}{\begin{corollary}}
 \newcommand{\EC}{\end{corollary}}
 \newcommand{\BR}{\begin{rem}}
 \newcommand{\ER}{\end{rem}}
 \newcommand{\BD}{\begin{defn}}
 \newcommand{\ED}{\end{defn}}
 \newcommand{\BI}{\begin{itemize}}
 \newcommand{\EI}{\end{itemize}}

\begin{document}
\title
[ Exponential mixing for some SPDEs with L\'evy noise] {Exponential mixing for some  SPDEs with L\'evy noise }
\author[E. Priola]{Enrico Priola}
\address{Dipartimento di Matematica, Universit\`a di Torino,
via Carlo Alberto 10 \\ 10123 Torino, Italy} \email{enrico.priola@unito.it}
\author[L.Xu]{Lihu Xu}
\address{TU Berlin, Fakult\"{a}t II, Institut f\"{u}r Mathematik,
Str$\alpha \beta$e des 17. Juni 136, D-10623 Berlin, Germany}
\email{xu@math.tu-berlin.de}
\author[J.Zabczyk]{Jerzy Zabczyk}
\address{Institute of Mathematics, Polish Academy of Sciences, P-00-950 Warszawa, Poland}
\email{zabczyk@impan.pl}

\thanks{ The first author gratefully acknowledges  the support  by the M.I.U.R. research project Prin 2008
``Deterministic and stochastic methods in the study of evolution
problems''. The second author gratefully acknowledges the support
from EURANDOM and Hausdorff Research Institute for Mathematics. The
third author gratefully acknowledges  the support by Polish Ministry
of Science and Higher Education grant ``Stochastic equations in
infinite dimensional spaces'' N N201 419039 }
%
\subjclass[2000]{}
\date{}
\maketitle
\begin{abstract} \label{abstract}
\noindent We show how gradient estimates for  transition semigroups
can be used to establish exponential mixing for a class of Markov
processes in infinite dimensions. We concentrate on semilinear
systems driven by cylindrical $\alpha$-stable noises  introduced in
\cite{PrZa09}, $\alpha \in (0,2)$. We  first  prove that if the
nonlinearity is bounded, then the system is ergodic and strong
mixing. Then we show that the system is   exponentially mixing
provided that the  nonlinearity, or its  Lipschitz constant, are
sufficiently small.
\end{abstract}

\section{Introduction}

\vskip 1mm The paper is concerned with asymptotic properties of the stochastic evolution equation:
\begin{equation}\label{e:MaiSPDE}
dX(t)=[AX(t)+F(X(t))]dt+dZ_t ,\,\,\, X_0=x,
\end{equation}
 driven by $\alpha$-stable, $\alpha \in (0,2)$, cylindrical process $Z$, introduced in \cite{PrZa09}. The paper \cite{PrZa09}
 investigated structural properties of $X$ like strong Feller and  irreducibility. Solutions turned out
 to be stochastically continuous but
 in general without c{\'a}dl{\'a}g modifications, see \cite{BGIPPZ}. The stochastic PDEs driven by L$\acute{e}$vy
 noises have been  intensively studied since some time, see e.g.  \cite{AWZ98},  \cite{AMR09}, \cite{PeZa06},
 \cite{Ok08}, \cite{MaRo09}, \cite{PrZa09} and \cite{XZ09}, and  the book \cite{PeZa07} for additional references.
 Even for equations like (\ref{e:MaiSPDE}), with the additive noise,  some basic questions are still open.

\vskip 1mm General results on existence of invariant measures for the linear equation
\begin{equation}
dX(t)=AX(t)dt+dZ_t ,\,\,\, X_0=x,
\end{equation}
with $A$ being the infinitesimal generator of a strongly continuous semigroup, were obtained  in \cite{Ch1987}, and more recently in \cite{PrZa09-2} and
\cite{FuXi09}, with different assumptions on $Z$. The nonlinear case was investigated in  \cite{Rusinek}, and by the dissipativity method in
\cite{PeZa07}; see also the recent paper \cite{XZ09}.

 In  the present paper we are concerned with the  exponential convergence to equilibrium,  called exponential mixing, see e.g.
 \cite[Proposition 16.5]{PeZa07} .We assume that the nonlinearity
  $F$ is Lipschitz continuous and bounded and  find explicit conditions on  $F$ under which the exponential mixing takes
 place. The speed of the convergence to equilibrium can be deduced from the proofs. Our result seems to be new even in the well-studied
 case  when $Z$ is a cylindrical Wiener process, corresponding to $\alpha =2$, requiring only obvious  modifications of our proofs.

 To prove the exponential mixing  we establish first gradient estimates for the transition semigroup of the process $X$ using the so called mild Kolmogorov equation. This version of the Kolmogorov equation was introduced in \cite{CDP} and elaborated in \cite{DPZ92} and \cite{DPZ96}. With a similar aim gradient estimates were used in \cite{XZ09}, in a different framework, with the assumption $\alpha \in (1,2)$ and without the exponential bounds. The ergodicity in \cite{XZ09}
and \cite{XZ09-1} was obtained by the technique of interacting
particle systems called finite speed of propagation of
information. For the finite dimensional stochastic system driven by
$\alpha$-stable processes, \cite{XZ09} obtained the strong mixing
result for small $F$.

\vskip 1mm Our main results are formulated as Theorems \ref{t:Erg}
and  \ref{t:ExpMix}. The first one establishes ergodicity without
additional  condition on $F$. Here it was enough to use the method
sketched in \cite{DPZ95}. Theorem \ref{t:ExpMix} is on exponential
mixing with additional restrictions on $F$. If  the nonlinearity is
only Lipschitz and bounded (and not small), it seems hard to obtain
exponential mixing by the  gradient estimates. We are now working on
this general case trying a different method.

\vskip 1mm The organization of the paper is as follows. Section
\ref{s:NM} introduces the notations and the main theorems. Section
\ref{s:PTE} proves the ergodicity and strong mixing for bounded
nonlinearity. The next section is devoted to exponential gradient
estimates. In the final section, we prove  exponential mixing for
small nonlinearities.

\section{Notations and main results} \label{s:NM}
We shall study our problem in a separable Hilbert space $H$ with an
orthonormal basis $\{e_k\}_{k \geq 1}$. Denote the inner product and
norm of $H$ by $\Ll\cdot,\cdot\Rr_H$ and $|\cdot|_H$ respectively.
 Let  $C_b(H,H)$ be the Banach space  of all bounded continuous functions
$f:H \rightarrow H$ with the supremum norm
$$||f||_{0}:=\sup_{x \in H} |f(x)|_H.
$$
Similarly, $C_b(H,\R)$ denotes the Banach space of all bounded continuous functions
$f:H \rightarrow \R$ with the supremum norm
$||f||_{0}:=\sup_{x \in H} |f(x)|.$

Given $f \in C_b(H,\R)$ and $x \in H$, for any $h \in H$, we write
\Be D_h f(x):=\lim_{\e \rightarrow 0} \frac{f(x+\e h)-f(x)}{\e} \Ee
provided the the above limit exists. If
 $h \mapsto D_h f(x)$ is linear and
$|D_h f(x)| \leq C|h|_H$, $h \in H$,
then there exists a unique element in $H$, denoted by $Df(x)$, such
that $D_h f(x)=\Ll D f(x),h \Rr_H$. Clearly, $Df$ is a function from
$H$ to $H$. Denote by $C^1_b(H,\R)$ the set of all bounded
differentiable functions $f: H \rightarrow \R$ with norm
$$||f||_{C^1_b}:=\sup_{x \in H} |f(x)|+\sup_{x \in H} |D f(x)|_H.$$
 Let $B_b(H,\R)$ be the set of all  bounded (Borel) measurable functions from $H$ into $\R$.

Further, we denote by $\mcl D$ the set of all bounded continuous
cylindrical functions $f: H \rightarrow \R$, i.e., $f$ is some
bounded continuous function depending on finite number of
coordinates of $H$. For any $f \in \mathcal{D}$, denote by
$\Lambda(f)$ the localization set of $f$, i.e. $\Lambda(f)$ is the
smallest set $\Lambda \subset \N$ such that $f \in
C_b(\R^{\Lambda},\R)$.
 Finally,  define $\mcl D^1:=\mcl D \cap C^1_b(H,\R)$.

 Let $z(t)$ be the normalized one dimensional symmetric $\alpha$-stable process, $0<\alpha<2$, with
 the following characteristic function
\Be \label{al} \E [e^{i \lambda
z(t)}]=e^{-t|\lambda|^{\alpha}},\,\,\lambda \in \R .
\Ee
The density of $z(1)$ will be denoted by $p_{\alpha}$. The infinitesimal generator $\partial^{\alpha}_x$  of the process is of the form
\begin{equation} \label{e:fraclap}
\partial^{\alpha}_x f(x)=
\frac{1}{C_{\alpha}} \int_{\mathbb{R} }
\frac{f(y+x)-f(x)}{|y|^{\alpha+1}}dy
\end{equation}
with  $C_{\alpha}= -\int_{\mathbb{R} }
(cosy-1)\frac{dy}{|y|^{1+\alpha}},$ see \cite{ARW00}.

\vskip 2mm We will consider equation (\ref{e:MaiSPDE}) under  the following assumptions:
\begin{assumption} \label{main}

 \BI

\ \quad

\item $A$ is a dissipative operator defined by
\begin{equation} \label{pu0}
A=\sum_{k \geq 1} (-\gamma_k) e_k \otimes e_k.
\end{equation}
$0<\gamma_1 \leq \gamma_2 \leq \cdots \leq \gamma_k \leq \cdots$ and
$\gamma_k \rightarrow \infty$ as $k \rightarrow \infty$.
\item
$F: H \rightarrow H$ is Lipschitz and bounded .
\item $Z= (Z_t)$ is a cylindrical $\alpha$-stable process with
$Z_t=\sum_{k \geq 1} \beta_k z_k(t) e_k$ where $\{z_k(t)\}_{k \geq
1}$ is a sequence of i.i.d. symmetric $\alpha$-stable processes (defined on some fixed
  stochastic basis
  $(\Omega, {\mathcal F}, ({\mathcal F}_t)_{t \ge 0}, {\PP}$)
 with
$0<\alpha<2$ and $\beta_k >0$, $k \ge 1,$ and
\begin{equation} \label{pu01}
\sum_{k \geq 1} \frac{\beta^{\alpha}_k}{\gamma_k}<+ \infty.
\end{equation}
\item There exists  $\sigma \in (0,1) $  such that
 \begin{equation} \label{pu1}
 \sup_{k \ge 1}{\frac{\gamma_k ^{{\frac{1}{\alpha}} -\sigma}}{\beta_k}} < +\infty.
\end{equation}
\EI
\end{assumption}

 \quad

\begin{remark}  Note that  $A$ generates a strongly continuous  semigroup of compact contractions $(e^{tA})$ on $H$.  Moreover, $e^{tA} e_k = e^{-t \gamma_k} e_k$, $k \ge 1$.

The first three conditions imply the existence of the  unique mild solution $X= (X(t,x))$ to
(\ref{e:MaiSPDE}), see \cite{PrZa09}. This is a predictable
   $H$-valued
   stochastic process, depending on  $x
 \in H$,  such that, for any $t \ge 0,$ $x \in H$,  it holds
($\PP$-a.s.): \begin{align} \label{mil}
 X(t,x) = e^{tA} x + \int_0^t e^{(t-s)A} F(X(s,x))ds
   + Z_A (t),\;\;\; \mbox{with}
   \;\; Z_A(t) = \int_0^t e^{(t-s)A}  d Z_s.
 \end{align}
 Condition \eqref{pu1}  is used
 in \cite[Section 5]{PrZa09} to  get  the strong Feller
 property for the  Markov semigroup $P_t$ associated with the solution $X$ of (\ref{e:MaiSPDE}):
$$
P_t f(x)=\E[f(X(t,x))],\,\,\,f\in B_b(H,\R), \;\;
 x \in H, \; t \ge 0.
$$
\end{remark}

For further use we need an equivalent formulation of (\ref{pu1}).

\BL \label{p:equivalence }  For arbitrary $\sigma >0$ the following two
conditions are equivalent:
\begin{equation}
B =
\sup_{k \ge 1}{\frac{\gamma_k ^{{\frac{1}{\alpha}} -\sigma}}{\beta_k}} < +\infty ,
\end{equation}
\begin{equation} \label{pu}
k_t =  \sup_{k \ge 1} \frac{ e^{- \gamma_k t }\, \gamma_k ^{1/ \alpha}} { \beta_k^{}}
 \le \frac{\hat c \, e^{ - \gamma_{1}  t /2  }} {t^{\sigma}}, \;\;
 t>0,
  \end{equation}
where\,\,\, $\hat c = B \, 2^{(2/\alpha) - \, \sigma }\, \frac{\sigma^{\sigma}}{e^{\sigma}}.$
\EL
\begin{proof}
To establish the equivalence  one can argue as in the beginning of \cite[Section 5]{PrZa09}.
 We only note the following estimate, for any $n \ge 1$,
$$
\frac{ e^{- \gamma_n t }\, \gamma_n ^{1/ \alpha}} { \beta_n^{}}
 \le 2^{{\frac {1}{\alpha}}} \exp (- \gamma_1 t /2 )
 \frac{ e^{- \frac{\gamma_n}{2} t }  \, {(\gamma_n/2) ^{1/ \alpha}}} {
\beta_n^{}}.
$$
\end{proof}

According to \cite{PrZa09}, we have the following result
for the system \eqref{e:MaiSPDE}.

\BT \label{p:MarSemPt} Under the  Assumption \ref{main} there exists a unique mild solution $X(t,x)$ for \eqref{e:MaiSPDE}. Moreover its associated transition semigroup $P_t$ is  strong Feller and irreducible.
 \ET
This paper aims to study the long time behaviour of the system
\eqref{e:MaiSPDE}.  The two main theorems are as follows.

\BT \label{t:Erg} Under Assumption \ref{main} there exists a unique invariant measure $\mu$ for the system \eqref{e:MaiSPDE}. The measure $\mu$ is strong
mixing, i.e.
$$\lim_{t \rightarrow \infty} P_t f(x)=\mu(f),$$
for all $f \in B_b(H,\R)$ and $x \in H$. \ET

In the formulation of the next theorem $p_{\alpha}$ stands for the density of $z(1)$ and the constant ${\hat c}$ was introduced in Lemma
\ref{p:equivalence }. \BT \label{t:ExpMix} Assume that Assumption \ref{main} holds and
 that one of the following two conditions holds:
\begin{enumerate}
\item [(i)]  $L_F<\gamma_1$,\,\,\,\,\,\,
where $L_F$ is the best  Lipschitz constant of $F$;
\\
\item [(ii)]
$||F||_{0}< {\frac{C_0}{\Gamma(1-\sigma)}} ({\frac{\gamma_1}{2}})^{1-\sigma}$,\,\, where $C_{0}= {\hat c}
\int_{\R}{\frac{(p^{'}_{\alpha}(z))^2}{p_{\alpha}(z)}} dz$ .
\end{enumerate}
Then the system \eqref{e:MaiSPDE} is exponentially mixing. More
precisely there exist constants $C=C(|x|_{H},\alpha,(\beta_n),
(\gamma_n), F)>0$ and $c=c(\alpha, (\beta_n),(\gamma_n),F)>0$ such
that  \Be \label{serve} \left|P_t f(x)-\mu(f)\right| \leq Ce^{-ct}
\, ||f||_{C^1_b},
 \Ee for all $f \in C^1_b(H, \R)$  and $x \in H$.
\ET
 We do not give explicit formulas for the constants $C$ and $c$ but they can be obtained by a careful examination of the proofs.

\smallskip Let us give some examples which the two main theorems can be applied to.

\begin{ex}
Consider the following stochastic semilinear equation on $D=[0,
\pi]^d$ with $d \geq 1$ with Dirichlet boundary conditions
\Be
\label{e:SemLin}
\begin{cases}
dX(t,\xi)=[\Delta X(t,\xi)+F(X(t,\xi))]dt+dZ_t(\xi), \\
X(0,\xi)=x(\xi), \\
X(t,\xi)=0, \ \ \xi \in \p D,
\end{cases}
\Ee
where $Z_t$ and $F$ are both specified below.
It is clear that $\Delta$ with Dirichlet boundary condition has the following
eigenfunctions
$$e_k(\xi)=\left(\frac{2}{\pi}\right)^{\frac d2} sin(k_1 \xi_1) \cdots sin(k_d \xi_d), \ \ \ \ k \in \N^d, \ \xi \in D.$$
It is easy to see that $\Delta e_k=-|k|^2 e_k$, i.e. $\gamma_k=|k|^2= k_1^2 + \ldots + k_d^2$,
  for all $k \in \N^d$.
We study the dynamics \eqref{e:SemLin} in the Hilbert space $H= L^2(D)$  with orthonormal basis $\{e_k\}_{k \in \N^d}$.
$Z= (Z_t)$ is some cylindrical $\alpha$-stable noises which, under the basis $\{e_k\}_k$, is
defined by
$$Z_t=\sum_{k \in \N^d} |k|^{\beta} z_k(t)e_k$$
where $\{z_k(t)\}_k$ are i.i.d.
 symmetric $\alpha$-stable processes
with $\alpha \in (0,2)$ and $\beta$ a real number. Note that
$\sum_{k \in \N^d} \frac{|k|^{\beta \alpha}}{|k|^2}< \infty$
if and
only if $2> d + {\alpha}{\beta}$.

\noindent From Theorems \ref{t:Erg}  and \ref{t:ExpMix}, we have

\BI
\item If $F$ is a bounded Lipschitz
function and
$$
2> d + {\alpha}{\beta}, \,\,\,   \frac{1}{\alpha}-
 \frac{\beta}{2}  < 1,
$$
or equivalently,
$$
{\frac{d}{\alpha}} < {\frac{2}{\alpha}} - \beta < 2,
$$
then the system \eqref{e:SemLin} is  strongly mixing. 
\item If in addition $||F||_{0}$ is sufficiently small then the system \eqref{e:SemLin} is  exponentially mixing.
\EI
 \end{ex}
\section{Proof of Theorem \ref{t:Erg}} \label{s:PTE}
According to Sections 5.2 and 5.3 of \cite{PrZa09}, the system \eqref{e:MaiSPDE} is irreducible and strong Feller. By Doob's theorem (see  \cite[Theorem 4.2.1]{DPZ96}), to prove Theorem \ref{t:Erg}, one only needs to show the existence of invariant measures. To this purpose  it is enough to establish
the tightness of $\{\mcl L(X(t,x))\}_{t\ge 1}$
 for some $x \in H$ (here
   $\mcl L(X(t,x))$ denotes the law of the random variable $X(t,x)$, see \eqref{mil}). For this we basically follow  \cite{DPZ95}.

Fix any $x \in H$ and set $X(t,x)= X(t)$. First note that for  any
fixed $t>0$, $e^{At}$ is a contraction  and a compact operator.
Contraction property is clear by the assumption of $A$. The
compactness is an easy corollary of the fact that the eigenvalues
$-\gamma_k$ tend to $-\infty$.

According to   \cite[Proposition 4.2]{PrZa09}, for any (small) $\e>0$, there exists some $M= M_{\epsilon}>0$
such that \Bes \PP\{|Z_A(t)|_H \leq M\} \geq 1-\e \Ees uniformly for all $t\ge 0$.  By the assumptions  on $F$ and $A$, one clearly has $|\int_0^t
e^{A(t-s)} F(X(s)) ds|_H \leq C$ uniformly for $t\ge 0$ and $\omega \in \Omega$. Hence,
\Be \label{e:XtBou}
 \PP\{|X(t)|_H \leq |x|+C+M\} \geq
\PP\{|Z_A(t)|_H \leq M\} \geq 1-\e.
\Ee
Let us rewrite \eqref{mil} as
 $$
 X(t)=e^{A} X(t-1)+\int_{t-1}^t e^{A(t-s)} F(X(s))ds+\int_{t-1}^t e^{A(t-s)} dZ_s,
  $$ thanks to the compactness of
$e^{A}$ and \eqref{e:XtBou}, the family  $\mcl L (\{e^{A} X(t-1)\})_{t\ge 1}$ is tight. The integrals $\int_{t-1}^t e^{A(t-s)} dZ_s$ have the same law as
$Z_A(1)$ for all $t\ge 1$, and thus their laws are  of course tight.

To complete the proof of  tightness of $\{\mcl L(X(t))\}_{t\ge 1}$, it is enough to show that the values of the integrals $\int_{t-1}^t e^{A(t-s)}
F(X(s))ds$ are contained in  a  compact set, for any $t \ge 1$. To this purpose note that the operator ${\mathcal R}$ from $L^{2}(0,1;H)$ into $H$
$$
{\mathcal R} \phi = \int_{0}^{1} e^{As} \phi (s)ds,\,\,\,\phi \in L^{2}(0,1;H),
$$
is compact (this follows from the compacteness of the operators
$e^{At},\,t>0$). Using the compactness of ${\mathcal R}$ and the
fact that  the transformation $F$ is bounded, we get that the
integrals are contained in a fixed compact set. The proof is
complete.

\section{Exponential gradient estimates}
To establish Theorem \ref{t:ExpMix} we  first derive exponential gradient estimates for the Galerkin approximation of the equation (\ref{e:MaiSPDE}).

\subsection{Galerkin approximation}
Let $\{e_k\}_{k \geq 1}$ be the orthonormal basis associated with the operator $A$ and write $\Gamma_N:=\{1,\ldots,N\}$

For any $x \in H$ and any integer $N>0$  we have the following
approximation of Eq. \eqref{e:MaiSPDE}:
\begin{equation} \label{e:GalApp}
\begin{cases}
d X^N_k(t)=[-\gamma_k X^N_k(t)+F^N_k(X^N(t))]dt+\beta_k dz_k(t), \\
X^N_k(0)=x_k,
\end{cases}
\end{equation}
for all $k \in \Gamma_N$, where $x^N=(x_k)_{k \in \Gamma_N}$, $x_k =
\langle x, e_k\rangle_H$, and $F^N_k(x^N)=\Ll F(x^N,0), e_k\Rr_H$.
Eq. \eqref{e:GalApp} can be written in the following vector form
\begin{equation} \label{e:AppEqu}
\begin{cases}
dX^N(t)=[A X^N(t)+F^N(X^N(t))]dt+dZ^N_t, \\
X^N(0)=x^N,
\end{cases}
\end{equation}
where $X^N(t)=(X^N_k(t))_{k \in \Gamma_N}$ and $Z^N_t=(\beta_k z_k(t))_{k \in \Gamma_N}$. The infinitesimal generator of \eqref{e:AppEqu} is
\Be
\label{e:AppEquGen}
\begin{split}
\mathcal{L}_N&=\sum_{k \in \Gamma_N} \beta^{\alpha}_k \p^{\alpha}_k+
\sum_{k \in \Gamma_N}[-\gamma_k x_k+F^N_k(x^N)] \p_k \\
&=\sum_{k \in \Gamma_N} [\beta^{\alpha}_k \p^{\alpha}_k-\gamma_k x_k \p_k]+\sum_{k \in \Gamma_N} F_k^N(x^N)\p_k,
\end{split}
\Ee where $\p_k=\p_{x_k}$ and $\p^{\alpha}_k=\p^{\alpha}_{x_k}$. In
the sequel we will identify  $x^N \in \R^N$ with
$$
 x^N = \sum_{k=1}^N x_k e_k \in H.
$$

Consider the Kolmogorov equation of the Galerkin approximation \Be \label{e:KolGalApp}
\begin{cases}
\p_t u^N(t)=\mcl {L}_N u^N(t), \\
u^N(0)=f,
\end{cases}
\Ee where $f \in \mcl D^1$ with $\Lambda(f) \subset \Gamma_N$. According to Section 5.3 of \cite{PrZa09}, Eq. \eqref{e:KolGalApp} has a mild solution
$P^N_t f$ which satisfies \Be \label{e:MilFor} P^N_t f(x^N) =
S^N_t f(x^N)+\int_0^t S^N_{t-s} [\Ll F^N, D P^N_s f\Rr_H](x^N)ds, \Ee
 where
$S^N_t$ is the Ornstein Uhlenbeck transition semigroup generated by the operator $\sum_{k \in \Gamma_N} [\beta^{\alpha}_k \p^{\alpha}_k-\gamma_k x_k
\p_k]$. Moreover, we also have
$$P^N_t f(x^N)=\E[f(X^N(t,x^N))].$$
>From the second step of the proof of   \cite[Theorem 5.7]{PrZa09}, we have
  \Be \label{e:LimSem} \lim_{N \rightarrow \infty} P^N_t f
(x^N)=P_tf(x), \Ee in particular,  for any $f \in
 \mcl D^1 $ with $\Lambda(f) \subset \Gamma_N$, $t \ge 0,$
 $x \in H$,
where $P_t$ is the transition semigroup defined in Theorem \ref{p:MarSemPt}.

\subsection{Estimates and their proofs} The gradient estimates which  are established here are of two different types. The first one is straightforward, and although
formulated  for (\ref{e:GalApp}),  is true in a  much more general
situation. Moreover,  it  is stated in terms of  the best Lipschitz
constant of the nonlinearity $F$. The second one is true for the
specific systems considered in the paper, requires more subtle
considerations and  is stated in terms of the supremum of $\{|F(x)|;
x\in H\}$. \BP \label{t:ExpDec0} Let $P^N_t$ be the transition
semigroup corresponding to the solution of the equation
(\ref{e:GalApp}). Then
\begin{equation} \label{e:EasDec}
||DP^N_t f||_{0}  \leq e^{-(\gamma_1-L_F)\, t} ||Df||_{0}, \ \
 \ f \in \mcl D^1 \;\;
  \text{with} \;
   \Lambda(f) \subset \Gamma_N,  \;\; t \ge 0.
\end{equation}
\EP \Bp
 Denoting by $X^N(t,x^N)$ the solution of Eq. \eqref{e:GalApp}
with initial data $x^N$, one  has
\begin{equation*} \label{e:TriCas}
|X^N(t,x^N)-X^N(t,y^N)|_H \leq e^{-(\gamma_1-L_F)t}|x^N-y^N|_H,
 \; t \ge 0, \; x,y \in H,
\end{equation*}
which implies $\lim_{\e \rightarrow 0}\sup |\frac{X^N(t,x^N +\e
h^N)-X^N(t,x^N)}{\e}| \leq e^{-(\gamma_1-L_F)t} |h^N|$, $h \in H$.
Hence, by the dominated convergence theorem one has
\begin{equation*}
\begin{split} |D_h P^N_t f(x^N)| & \leq \E\left[|Df(X^N(t))|_H \lim_{\e
\rightarrow 0}\sup \left|\frac{X^{N}(t,x^N+\e
h^N)-X^N(t,x^N)}{\e}\right|\right] \\
&  \leq e^{-(\gamma_1-L_F)\, t} ||Df||_{0} \, |h|_H.
\end{split}
\end{equation*}
\Ep

\BP  \label{t:ExpDec} Let $P^N_t$ be the transition semigroup
corresponding to the solution of the equation (\ref{e:GalApp}).
There exists a positive constant $C$,\,depending on $\alpha,
(\beta_n), (\gamma_n)$, $\sigma$, such that \Bes ||DP^N_tf||_{0}
\leq C e^{-\omega t} ||Df||_{0}, \ \  \ f \in \mcl D^1, \;\;
\text{with} \; \;  \Lambda(f) \subset \Gamma_N, \; \; t \ge 0, \Ees
where
$$
\omega = {\frac{\gamma_1}{2}} - (C_{0}\| F\|_{0}\Gamma(1 - \sigma))^{{\frac{1}{1-\sigma}}},
$$
and the constant $C_{0}$ was introduced in Theorem \ref {t:ExpMix}. \EP

To prove Proposition \ref{t:ExpDec}, we need the following lemma. \BL
 \label{l:DecSNt}
  For any $f \in \mcl D^1$, with $\Lambda(f) \subset \Gamma_N,$ one has:
 \begin{align*}
& (i) \; |DS^N_tf(x^N)|_H \leq  e^{-\gamma_1 t} ||Df||_{0}, \\
& (ii)\; |DS^N_tf(x^N)|_H \leq  \frac{C_0 e^{ - \gamma_{1}  t /2  }} {t^{\sigma}} ||f||_{0},
\end{align*}
for all $x \in H$, $N \ge 1$, $t>0$
 (where $C_0 = {\hat c} \int_{\R}{\frac{(p^{'}_{\alpha}(z))^2}{p_{\alpha}(z)}} dz$, see Theorem \ref{t:ExpMix}).
  \EL \Bp
 Note that, for any vector $h \in \sum_{k=1}^N h_k e_k$,
  $x \in H,$ we have
$$
 \langle DS^N_tf(x^N), h \rangle =
 \langle DS_t f(x^N), h \rangle,
$$
 where $S_t$ is the Ornstein-Uhlenbeck semigroup acting on $C_b (H,\R)$
  associated to $X(t,x)$ when $F=0$. We know that
 \begin{equation} \label{ou}
 S_t f(x) =
    \int_H f(e^{tA} x +  y ) \, \,
  \mu^0_t (dy),
 \end{equation} where
 $\mu^0_t$ is the law of the random variable
  $Z_A(t) =  \int_0^t e^{(t-s)A} d Z_s$. By  differentiating under the
  integral sign in \eqref{ou} we immediately get the first assertion
 (remark that $ \| e^{tA} \| \le e^{- \gamma_1 t}$, $t \ge 0$, where $\| e^{tA} \|$ denotes the operator norm of $ e^{tA}$).

\smallskip As for the second assertion, recall the gradient estimate from \cite[Theorem 4.14]{PrZa09}
$$
|DS_tf(x)|_H \leq  8 c_{\alpha}  \Big( \sup_{k \ge 1} \frac{ e^{- \gamma_k t }\, \gamma_k ^{1/ \alpha}} { \beta_k^{}} \Big) ||f||_{0}, \; \,\,f \in C_b
(H, \R),\,\,\,\,t>0,
$$
where
$$
c_{\alpha} = {\frac{1}{8}}\int_{\R}{\frac{(p^{'}_{\alpha}(z))^2}{p_{\alpha}(z)}} dz.
$$
According to \eqref{pu} we have
\begin{align} \label{grad}
|DS_tf(x)|_H \leq  8 c_{\alpha} \frac{\hat c \, e^{ - \gamma_{1}  t /2  }} {t^{\sigma}}   ||f||_{0}, \; t>0,
\end{align}
and the assertion follows.
 \Ep
\Bp [Proof of Proposition \ref{t:ExpDec}] By \eqref{e:MilFor} and Lemma \ref{l:DecSNt}, we have, for any $N \in \N$,   \Bes
\begin{split}
|DP^N_t f(x^N)|  \leq e^{-{\frac{\gamma_1} {2}} t}||Df||_{0} +\int_0^t \frac{C_0 e^{-{\frac{\gamma_1} {2}} (t-s)}}{(t-s)^{\sigma}} ||\Ll F^N, DP^N_s
f\Rr_H||_{0} ds,
\end{split}
\Ees where $C_0$ is defined in Lemma \ref{l:DecSNt}; therefore,  \Bes
\begin{split}
||DP^N_t f||_{0}  \leq e^{-{\frac{\gamma_1} {2}} t}||Df||_{0} +\int_0^t \frac{C_0 e^{-{\frac{\gamma_1} {2}} (t-s)}}{(t-s)^{\sigma} } ||F||_{0} ||DP^N_s
f||_{0} ds,
\end{split}
\Ees Writing $v^N(t)=e^{{\frac{\gamma_1} {2}} t} ||DP^N_t f||_{0}$,
 we have from the above inequality \Bes v^N(t) \leq
||Df||_{0}+\int_0^t \frac{C_0 ||F||_{0}}{(t-s)^{\sigma}} v^N(s) ds. \Ees Now we use the following Henry's estimate (see
 \cite{henry}):
 let  $a \ge 0$, $b \ge 0$ and $\beta>0$ and consider
  a   non-negative
 locally integrable function $u$ on $[0,T)$ such that
$$
 u(t) \le a + b  \int_0^t (t-s)^{\beta -1} u(s) ds,
 \;\;\; t \in [0,T),
$$
then we have
 $u(t) \le a G_{\beta }(\theta t)$, $t \in [0,T)$, where
$$
 \theta = (b \Gamma(\beta))^{1 /\beta},\;\;\;
  G_{\beta}(z) = \sum_{n \ge 0} \frac{z^{n \beta}}{
  \Gamma (n \beta +1)}, \;\; z \ge 0
$$
(note that $G_{\beta } (z) \sim \frac{1}{\beta} e^{z}$
 as $z \to +\infty$).

In our case $u(t) = v^N(t)$, $a = \| Df\|_{0}$, $b =
  C_0 \| F\|_{0}$ and $\beta  = 1- \sigma$. Thus for a constant $C_1$, depending on $\sigma$, and all positive $t$
 we get
$$
v^N(t)  \le {\frac{C_1 \| Df\|_{0} }{ 1 - \sigma}} \; \exp \left\{\left(C_0
 \| F\|_{0} \Gamma(1 - \sigma)\right)^{ \frac{1}{1-\sigma}}\;  t
 \right\}
$$
Therefore
$$
||DP^N_t f||_{0} \le  {\frac{C_1 \| Df\|_{0}}{ 1 - \sigma}} \exp \left\{\left( C_0
 \| F\|_{0} \Gamma(1 - \sigma)\right)^{ \frac{1}{1-\sigma}}\; t
 -\frac{\gamma_1}{2}t\right \}
$$
and  we get the assertion. \Ep

\section{Proof of Theorem \ref{t:ExpMix}} \label{s:PTEM}
 From Theorem \ref{t:Erg},  the system (\ref{e:MaiSPDE})
is ergodic
 and has a unique invariant measure $\mu$.

Note that it is enough to prove
   \eqref{serve} for any $f \in  \mcl D^1$. Indeed then, approximating any function $f \in C^1_b (H)$ by a sequence $(f_n ) \subset {\mcl D^1}$ such that
  $f_n \to f$ and $Df_n \to Df$ pointwise with
   $\sup_{n \ge 1}\| f_n\|_{C^1_b} < \infty$, we
 get easily the complete assertion.
 Let us fix $f \in \mcl D^1 $ and suppose that
  $\Lambda(f) \subset \Gamma_N$.

The crucial point  of the proof is to apply the gradient estimate from Proposition \ref{t:ExpDec}, in the spirit of  \cite[Proposition 16.4]{PeZa07}.

Concerning the case (i), by using  \eqref{e:EasDec}, we can show the exponential mixing by a similar argument as for the case (ii).
So from now on we shall concentrate on the proof of the case (ii).

\smallskip
 First let us suppose that $\alpha \in (1,2)$. Let
$$
Z^N_A(t)=\int_0^t e^{A(t-s)} dZ^N_s
$$
 and write $Y^N(t)=X^N(t, x^N)-Z^N_A(t)$ so that
$$Y^N(t)=e^{At}x^N+\int_0^t e^{A(t-s)} F^N(X^N(s, x^N))ds.
$$
For any $t_2>t_1>0$,
writing $s=t_2-t_1$, by the gradient estimates of  Proposition
\ref{t:ExpDec}, we have, for any $x \in H$,

\Be \label{e:PNt20-PNt10}
\begin{split}
|P^N_{t_2} f(x^N)-P^N_{t_1} & f(x^N)|= | \E[P^N_{t_1}f(X^N(s,x^N))
-P^N_{t_1} f(x^N)]| \\
&\leq ||DP^N_{t_1}f||_{0}\,
   \E | X^N(s,x^N) - x^N |_H \\
&\leq Ce^{-\omega t_1} ||Df||_{0}\, \E | X^N(s,x^N) - x^N |_H,
\end{split}
\Ee  where $C,\omega>0$ are given in Proposition
 \ref{t:ExpDec}.  Now note that
$$
| X^N(t,x^N) - x^N |_H \le | e^{A t} x^N- x^N |_H +
 \big| \int_0^t
e^{A(t-s)} F^N(X^N(s, x^N)ds \big |_H + |Z_A(t)|_H.
$$
We have
$$
\big| \int_0^t e^{A(t-s)} F^N(X^N(s, x^N))ds \big |_H
 \le \| F \|_{0} \int_0^{\infty}  e^{ - \gamma_1 s} ds
  \le \frac{\| F \|_{0} }{\gamma_1},
$$
 for any $\omega \in \Omega$, $x \in H$ and $t \ge 0$. Concerning
 $|Z_A(t)|_H$ remark that since $\alpha >1$
 (see  \cite[Theorem 4.4]{PrZa09})
$$
 \E |Z_A(t) |_{H}  \le  \tilde c_1  \Big (  \sum_{n \ge
1} |\beta_n|^{\alpha} \frac{(1- e^{- \alpha \gamma_n t })}{\alpha
 \gamma_n }
\Big)^{1/\alpha} \le \tilde c_1  \Big (  \sum_{n \ge 1}
|\beta_n|^{\alpha} \frac{1}{\alpha
 \gamma_n }
\Big)^{1/\alpha} < +\infty, \;\; t \ge 0.
$$
 It follows that
$$
\E | X^N(s,x^N) - x^N |_H \le C_2 (1 + |x^N|_H), \;\;
$$
 where $C_2$ does not depend on $N$, $s \ge 0 $ and $ x \in H$.

\vskip 2mm
 By  \eqref{e:PNt20-PNt10}
  we get
  $$
|P^N_{t_2} f(x^N)-P^N_{t_1}  f(x^N)|
 \leq C_3 e^{-\omega t_1} ||Df||_{0} \,  (1+ |x^N |_H).
$$
 Passing to the limit as $N \to \infty$, we get (see the proof
  of Theorem 5.7 in \cite{PrZa09})
$$
|P_{t_2} f(x)-P_{t_1}  f(x)|
 \leq C_3 e^{-\omega t_1} ||Df||_{0} \,  (1+ |x |_H).
$$
where $C_3$ does not depend on $t_1$. This estimate shows that
 $P_tf(x)$  converges to some constant  exponentially fast
  as $t \to +\infty$. By
ergodicity of the system, this constant must be $\mu(f)$. \vskip 2mm
Let us consider now $\alpha \in (0,1]$. For any $f \in C^1_b(H,\R)$, $s \in (0,1]$,
define $[f]_{s}$ by
$$
 [ f]_{s}:=\sup_{x\neq y }
 \frac{|f(x)-f(y)|}{|x-y|^{s}}.
$$
We have
\begin{equation} \label{kll}
[ f]_{s} \leq 2^{1-s} \sup_{x\neq y }
 \frac{|f(x)-f(y)|^s}{|x-y|^{s}} ||f||^{1-s}_{0} \leq 2^{1-s} ||f||^{1-s}_{0} ||Df||^s_{0}.
\end{equation}
Now we choose  $p \in (0, \alpha)$.  Using
that  $\E |Z_A(t) |_{H}^p < \infty $ (see \cite[Theorem 4.4]{PrZa09})
 and \eqref{kll} with $s=p$ we get, arguing as before,
\begin{align*}
|P^N_{t_2} f(x^N)-P^N_{t_1} & f(x^N)|=| \E[P^N_{t_1}f(X^N(s,x^N))
-P^N_{t_1} f(x^N)]| \\
&\leq 2^{1-p} \, ||DP^N_{t_1}f||_{0}^{p} \, \| f\|^{1-p}_{0}
   \; \E | X^N(s,x^N) - x |_H^p \\
&\leq 2^{1-p} \, C^p e^{- \omega p \, t_1} ||f||_{C^1_b}\;
 \E | X^N(s,x^N) - x^N |_H^p.
\end{align*}
    Since
$$
 \E |Z_A(t) |_{H}^p  \le  \tilde c_p  \Big (  \sum_{n \ge
1} |\beta_n|^{\alpha} \frac{(1- e^{- \alpha \gamma_n t })}{\alpha
 \gamma_n }
\Big)^{p/\alpha} \le \tilde c_p  \Big (  \sum_{n \ge 1}
|\beta_n|^{\alpha} \frac{1}{\alpha
 \gamma_n }
\Big)^{p/\alpha} < +\infty, \;\; t \ge 0,
$$
 it follows that
$$
\E | X^N(s,x) - x^N |_H^p \le C_4 (1 + |x^N|_H)^p, \;\;
$$
 where $C_4$ does not depend on $N$, $s \ge 0 $ and $ x \in H$. Finally we have
$$
|P_{t_2} f(x)-P_{t_1}  f(x)|
 \leq C_5 e^{-\omega p \,  t_1} ||f||_{C^1_b} \,  (1+ |x |_H)^p,
$$
where $C_5$ does not depend on $t_1$. Arguing as before we complete the proof.



\bibliographystyle{amsplain}
\

 \vskip  3mm \noindent   \textit{Acknowledgements.}
  The authors thank  the Newton Institute (Cambridge), where this paper was initiated,  for its hospitality.

\end{document}